\documentclass[a4paper,11pt]{article}
\AtBeginDvi{}
\usepackage{amsmath,graphicx}

\renewcommand{\title}[1]{\vspace{\fill}
\eject\addtolength{\baselineskip}{4pt}
{\bfseries\LARGE #1}\\[3mm]\addtolength{\baselineskip}{-4pt}}
\renewcommand{\author}[3]{\parbox[t]{75mm}
{\begin{center}{\scshape #1}\\[3mm] #2\\
 {\ttfamily #3} \end{center}}}

\usepackage{amssymb}
\usepackage{amsthm}
\usepackage{bm}
\usepackage{enumerate}
\usepackage[monochrome]{color}
\usepackage{mathtools}

\makeatletter
\newtheoremstyle{myplain}{\topsep}{\topsep}{\itshape}{}{\bfseries}{}{1em}
  {\thmname{\@ifempty{#2}{\@ifempty{#3}{#1}{#3}}{#1}}%
   \thmnumber{\@ifnotempty{#1}{ }\@upn{#2}}%
   \thmnote{\@ifnotempty{#2}{ {\the\thm@notefont(#3)}}}}
\makeatother
\theoremstyle{myplain}
\newtheorem{thm}{Theorem}
\newtheorem{lem}[thm]{Lemma}
\newtheorem{prop}[thm]{Proposition}

\newtheorem{cl}[thm]{Claim}

\makeatletter
\renewenvironment{proof}[1][\proofname]{\par
  \pushQED{\qed}%
  \normalfont \topsep6\p@\@plus6\p@\relax
  \trivlist
  \item[\hskip\labelsep
        \scshape
    #1\@addpunct{:}]\ignorespaces
}{%
  \popQED\endtrivlist\@endpefalse
}
\renewcommand{\qed}{%
  \ifmmode \mathqed
  \else
    \leavevmode\unskip\penalty9999 \hbox{}\nobreak
    \quad\hbox{\qedsymbol}%
  \fi
}
\renewcommand{\qedsymbol}{\textsquare}
\makeatother

\newcommand{\vct}{\bm}

\newcommand{\RR}{\mathbb{R}}
\newcommand{\ZZ}{\mathbb{Z}}

\newcommand{\CutP}{\mathrm{CUT}^\square}

\newcommand{\CorP}{\mathrm{COR}^\square}

\newcommand{\Gr}{\mathrm{Gr}}

\newcommand{\trans}{\mathrm{T}}
\newcommand{\K}{\mathrm{K}}
\newcommand{\A}{\mathrm{A}}
\newcommand{\B}{\mathrm{B}}
\newcommand{\Z}{\mathrm{Z}}
\newcommand{\I}{\mathrm{I}}
\newcommand{\N}{\mathrm{N}}

\newcommand{\symdiff}{\mathbin{\triangle}}

\newcommand{\revddots}
  {\mathinner{\mkern1mu\raise1pt
   \vbox{\kern7pt\hbox{.}}\mkern2mu
   \raise4pt\hbox{.}\mkern2mu\raise7pt\hbox{.}\mkern1mu}}

\begin{document}
\newlength{\parindentorig}
\setlength{\parindentorig}{\parindent}
\begin{center}

\title{New classes of facets of cut polytope \\
  and tightness of $\I_{mm22}$ Bell inequalities}
\author{David Avis}{
  School of Computer Science \\
  McGill University \\
  3480 University St. \\
  Montreal, Quebec \\
  Canada H3A 2A7
}{avis@cs.mcgill.ca}
\author{\underline{Tsuyoshi Ito}}{
  Dept.\ of Computer Science \\
  Grad.\ School of Info.\ Sci.\ and Tech. \\
  The University of Tokyo \\
  7-3-1 Hongo,
  Bunkyo-ku, Tokyo \\
  113-0033 Japan
}{tsuyoshi@is.s.u-tokyo.ac.jp}

\end{center}

\begin{quote}
  \textbf{Abstract:}
    The Grishukhin inequality $\Gr_7$ is a facet of $\CutP_7$, the cut
    polytope on
    seven points, which is ``sporadic'' in the sense that its proper
    generalization has not been known.
    In this paper, we extend $\Gr_7$ to an inequality $I(G,H)$ valid
    for $\CutP_{n+1}$ where $G$ and $H$ are graphs with $n$ nodes
    satisfying certain conditions, and prove a necessary and
    sufficient condition for $I(G,H)$ to be a facet.
    This result combined with the triangular elimination theorem of Avis,
    Imai, Ito and Sasaki settles Collins and Gisin's conjecture in
    quantum theory affirmatively: the $\I_{mm22}$ Bell
    inequality is a facet of the correlation
    polytope $\CorP(\K_{m,m})$ of the complete bipartite graph
    $\K_{m,m}$ for all $m \ge 1$.
    We also extend the $\Gr_8$ facet inequality of $\CutP_8$ to an
    inequality $I'(G,H,C)$ valid for $\CutP_{n+2}$, and provide a
    sufficient condition for $I'(G,H,C)$ to be a facet.
\end{quote}

\begin{quote}
\end{quote}
\vspace{5mm}

\section{Introduction}

Cut polytopes are convex polytopes which arise in many different
fields~%
\cite{DezLau-JCAM94:applications1,DezLau-JCAM94:applications2,%
DezLau:cut97}.
Since testing membership in cut polytopes is
NP-complete~\cite{AviDez-Net91}, it is unlikely that there exists a
concise and complete description of their facial structure in general.
Much efforts has been devoted to identifying classes of inequalities
which are valid for cut polytopes and have good properties.
Hypermetric, clique-web and parachute inequalities are examples of
classes of valid inequalities 
for which important subclasses are facet inducing.
For $N\le6$, all facets of $\CutP_N$, the cut polytope of complete
graph $\K_N$, are hypermetric.
However, $\CutP_7$ has a facet called the \emph{Grishukhin inequality}
$\Gr_7$ which is not known to belong to any such general class.
Efforts have been made to relate $\Gr_7$ to other
inequalities.
As a result, De~Simone, Deza and Laurent~\cite{DesDezLau-DM94} showed
that $\Gr_7$ is a collapse of a pure facet inequality $\Gr_8$ of
$\CutP_8$.

The cut polytope $\CutP(\K_{1,m,m})$ of the complete tripartite graph
$\K_{1,m,m}$ is linearly isomorphic to the correlation polytope
$\CorP(\K_{m,m})$ of the complete bipartite graph $\K_{m,m}$.
In quantum theory, the correlation polytope
$\CorP(\K_{m,m})$ is seen as the set of possible results of a series of Bell
experiments with a non-entangled (separable) quantum state shared by
two distant parties, where each party has $m$ choices of measurements.
In this context, a valid inequality of $\CorP(\K_{m,m})$ is called
a \emph{Bell inequality} and if facet inducing,
a \emph{tight Bell inequality}.
Readers are referred to \cite{WerWol-QIC01} for further information
about Bell inequalities.
Collins and Gisin~\cite{ColGis-JPA04} found a class of $\I_{mm22}$
inequalities valid for $\CorP(\K_{m,m})$ for general $m$ and
conjectured that for all $m \ge 1$, $\I_{mm22}$ inequality is tight, or
equivalently, that it is a facet of $\CorP(\K_{m,m})$.

Avis, Imai, Ito and Sasaki~\cite{AviImaItoSas:0404014} introduced an
operation called \emph{triangular elimination} to convert a facet of
$\CutP_N$ to a facet of $\CutP(\K_{1,m,m})$ for appropriate $m$.
By using this operation, the tightness of the $\I_{3322}$ and
$\I_{4422}$ Bell inequalities follows from the fact that the pure
pentagonal and the Grishukhin inequalities are facets of $\CutP_5$ and
$\CutP_7$, respectively.
This suggests that some natural extensions of the pure pentagonal and
the Grishukhin inequalities may give facets of $\CutP_{2m-1}$ for
$m\ge3$.
We will prove that it is the case and that hence the conjecture by
Collins and Gisin is true.
More specifically, we will introduce inequalities $I(G,H)$ valid for
$\CutP_{n+1}$ where $G$ and $H$ are graphs with $n$ nodes which satisfy
certain conditions described later, and prove a necessary and
sufficient condition for $I(G,H)$ to be a facet.

As further extensions, we apply to $I(G,H)$ an operation similar to
the one used to construct $\Gr_8$ from $\Gr_7$.
Actually this operation gives inequalities $I'(G,H,C)$ valid for
$\CutP_{n+2}$ where $C$ is a cycle of length four in $G$.
We will give a sufficient condition for $I'(G,H,C)$ to be a facet,
generalizing the fact that $\Gr_8$ is a facet of $\CutP_8$.

The rest of the paper is organized as follows.
In Section~\ref{sect:preliminaries}, we review the tools used later.
In Section~\ref{sect:ineq-1}, we introduce the inequality $I(G,H)$
valid for the cut polytope, which is a generalization of the $\Gr_7$
inequality, and we prove a necessary and sufficient condition for it
to be a facet.
Section~\ref{sect:ineq-2} defines the valid inequality
$I'(G,H,C)$, which is a generalization of the $\Gr_8$ inequality,
and we provide a sufficient condition for it to be a facet.
The proof of the sufficient condition is deferred to appendix.
In Section~\ref{sect:imm22}, we prove the tightness of $\I_{mm22}$
Bell inequalities.

\section{Preliminaries} \label{sect:preliminaries}

\subsection{Cut polytopes}

Here we review the definition of and results on cut polytopes only
briefly.
Readers are referred to the book by Deza and
Laurent~\cite{DezLau:cut97} for details.

\paragraph*{Definition}
The \emph{cut polytope} $\CutP(G)$ of a graph $G=(V,E)$ is a convex
polytope in
the vector space $\RR^E$ defined as the convex hull of the $2^{\lvert
V\rvert-1}$ different cut vectors $\vct{\delta}_G(S)$ for $S\subseteq V$.
The cut vector $\vct{\delta}_G(S)\in\RR^E$ is a 0/1 vector defined by
$\delta_{uv}(S)=1$ if and only if exactly one of $u$ and $v$ belongs
to $S$, where $uv$ denotes the edge connecting two nodes $u$ and $v$.
The cut polytope $\CutP(\K_N)$ of the complete graph $\K_N$ is denoted
by $\CutP_N$.

Similarly, the \emph{correlation polytope} $\CorP(G)$ is a convex
polytope in
$\RR^{V\cup E}$ defined as the convex hull of the $2^{\lvert V\rvert}$
correlation vectors $\vct{p}_G(S)$ for $S\subseteq V$.
The correlation vector $\vct{p}_G(S)\in\RR^{V\cup E}$ is a 0/1 vector
defined by $p_u(S)=1$ if and only if $u\in S$ and $p_{uv}(S)=1$ if and
only if $\{u,v\}\subseteq S$.

The correlation polytope $\CorP(G)$ of a graph $G=(V,E)$ is linearly
isomorphic to $\CutP(\nabla G)$, where $\nabla G$ is the
\emph{suspension graph} of $G$: the
graph obtained by adding to $G$ a new node $\Z$ adjacent to all the
nodes of $G$.
The linear isomorphism between them is called the \emph{covariance
mapping}:
$p_u=x_{\Z u}$ for $u\in V$ and $p_{uv}=\frac12(x_{\Z u}+x_{\Z v}
-x_{uv})$ for $u,v\in V$, $u\ne v$.

\paragraph*{Hypermetric inequalities}
Let $N\ge3$ be an integer and $\vct{b}\in\ZZ^N$ an integer vector with
$\sum_{i=1}^N b_i=1$.
The inequality
$\sum_{1\le i<j\le N} b_ib_jx_{ij}\le0$
is valid for $\CutP_N$ and called the \emph{hypermetric inequality}
defined by the vector $\vct{b}$.

While an exact characterization of when a hypermetric inequality
becomes a facet of $\CutP_N$ is not known, many sufficient conditions
are known.
We review here some of them which we use later.

\begin{thm}[Corollary~28.2.5~(i) in \cite{DezLau:cut97}]
    \label{thm:hyp-pure}
  Let $s\ge1$ be an integer, and $\vct{b}\in\ZZ^N$ be an integer
  vector with $s+1$ entries equal to $1$, $s$ entries equal to $-1$
  and the other $N-(2s+1)$ entries equal to $0$.
  Then the hypermetric inequality defined by $\vct{b}$ is a facet of
  $\CutP_N$.
  This inequality is called a \emph{pure $(2s+1)$-gonal inequality},
  or if $s=1$, simply a \emph{triangle inequality}.
\end{thm}

We define $T(u,v;w)=x_{uv}-x_{uw}-x_{vw}$.
By using this notation, a triangle inequality is written as
$T(u,v;w)\le0$.

\begin{thm}[``If'' part of Theorem~28.2.4~(iiib) in \cite{DezLau:cut97}]
    \label{thm:hyp}
  The hypermetric inequality defined by $\vct{b}$ with
  $b_1=\dots=b_{N-2}=1$, $b_{N-1}=-1$ and $b_N=-N+4$ is a facet of
  $\CutP_N$.
\end{thm}

\paragraph*{Switching of inequality}
We mention three operations on inequalities valid for cut polytopes.
One is the switching operation.
Let $G=(V,E)$ be a graph, $\vct{a}\in\RR^E$ and $a_0\in\RR$.
The \emph{switching} of the inequality $\vct{a}^\trans\vct{x}\le a_0$
by the cut $S\subseteq V$ is an inequality $\vct{b}^\trans\vct{x}\le
b_0$ with $b_{ij}=(-1)^{\delta_{ij}(S)}\cdot a_{ij}$ and
$b_0=a_0-\vct{a}^\trans\vct{\delta}_G(S)$.

Switching is an automorphism of the cut polytope $\CutP(G)$.
Therefore $\vct{b}^\trans\vct{x}\le b_0$ is valid (resp.\ a facet) if
and only if $\vct{a}^\trans\vct{x}\le a_0$ is valid (resp.\ a facet).

\paragraph*{Collapsing and lifting of inequality}
The other two operations are collapsing and lifting.
Let $G=(V,E)$ be a complete graph on node set $V$ and $uv\in E$.
Let $G'=(V',E')$ be the complete graph on node set
$V'=(V\setminus\{u,v\})\cup\{w\}$ with a new node $w$.

The \emph{$(u,v)$-collapse} of a vector $\vct{a}\in\RR^E$ is a vector
$\vct{a}^{u,v}\in\RR^{E'}$ defined by
\begin{alignat*}{2}
  a^{u,v}_{ij} &= a_{ij} &
    \quad &\text{for $i,j\in V\setminus\{u,v\}$, $i\ne j$,} \\
  a^{u,v}_{wi} &= a_{ui}+a_{vi} &
    \quad &\text{for $i\in V\setminus\{u,v\}$.}
\end{alignat*}
For $\vct{a}\in\RR^E$ and $a_0\in\RR$, an inequality
$(\vct{a}^{u,v})^\trans\vct{x}\le a_0$ is said to be the
\emph{$(u,v)$-collapse} of the inequality $\vct{a}^\trans\vct{x}\le
a_0$.

If the inequality $\vct{a}^\trans\vct{x}\le a_0$ is valid for
$\CutP(G)$, its collapse $(\vct{a}^{u,v})^\trans\vct{x}\le a_0$ is
valid for $\CutP(G')$.

The opposite operation of collapsing is called \emph{lifting}.
The following lemma provides a sufficient condition for lifting to
preserve a facet.
The proof of the lemma is given below Lemma~26.5.3 in the
book~\cite{DezLau:cut97}.

\begin{lem}[Lifting lemma~\cite{DezLau:cut97}]
    \label{lem:lifting}
  Let $\vct{a}\in\RR^E$.
  The inequality $\vct{a}^\trans\vct{x}\le0$ is a facet of $\CutP(G)$
  if the following conditions are satisfied.
  \begin{enumerate}[(i)]
    \item
      The inequality $\vct{a}^\trans\vct{x}\le0$ is valid for
      $\CutP(G)$, and its $(u,v)$-collapse
      $(\vct{a}^{u,v})^\trans\vct{x}\le0$ is a facet of $\CutP(G')$.
    \item
      There exist $\lvert V\rvert-1$ subsets $T_j$ of $V$ with
      $u\notin T_j$ and $v\in T_j$ such that the cut vectors
      $\vct{\delta}_G(T_j)$ are roots (vertices lying on the face) of
      $\vct{a}^\trans\vct{x}\le0$
      and the incidence vectors of $T_j$ are linearly independent.
  \end{enumerate}
\end{lem}

\begin{figure}
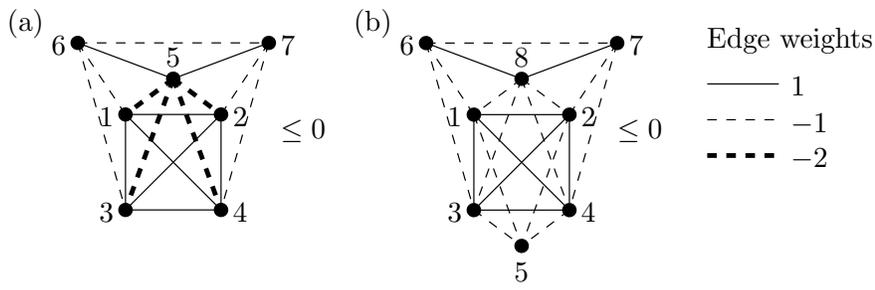

  \centering
  \begin{tabular}{ccc}
    (a) \raisebox{-\height}{\input{figures/gr7.pstex_t}}%
        \hspace*{1em}
    &
    (b) \raisebox{-\height}{\input{figures/gr8.pstex_t}}%
        \hspace*{1em}
    &
        \raisebox{-\height}{\input{figures/edge-weights-1m1m2.pstex_t}}%
        \hspace*{1em}
  \end{tabular}
  \caption{(a)~The Grishukhin inequality $\Gr_7$, which is a facet of
    $\CutP_7$.
    (b)~The $\Gr_8$ inequality, which is a facet of $\CutP_8$.}
  \label{fig:grishukhin}
\end{figure}

\paragraph*{Grishukhin inequality}
The cut polytope $\CutP_7$ has 11 inequivalent facets under permutation
and switching symmetries~\cite{Gri-EJC90,DesDezLau-DM94}.
All but one of them belong to at least one of three general classes of
valid inequalities: hypermetric, clique-web and parachute
inequalities.
The remaining facet is not known to belong to any classes that are as
general as these classes.
This ``sporadic'' facet is called the \emph{Grishukhin inequality}
$\Gr_7$.
The Grishukhin inequality looks like
$
  \sum_{1\le i<j\le4} x_{ij}
  +x_{56}+x_{57}-x_{67}-x_{16}-x_{36}-x_{27}-x_{47}
  -2\sum_{1\le i\le4} x_{5i} \le 0
$
and illustrated in Figure~\ref{fig:grishukhin}~(a).

De~Simone, Deza and Laurent~\cite{DesDezLau-DM94} found a facet of
$\CutP_8$ which is pure (all the coefficients are $0$ or $\pm1$) and
is a lifting of $\Gr_7$.
This facet is called $\Gr_8$ in \cite{DezLau:cut97} and illustrated in
Figure~\ref{fig:grishukhin}~(b).

\subsection{Bell inequalities}

\paragraph*{$\I_{mm22}$ Bell inequalities}
Collins and Gisin~\cite{ColGis-JPA04} showed that the $\I_{mm22}$
inequalities:
\begin{equation}
  -p_{\A_1}-\smashoperator{\sum_{1\le j\le m}} (m-j)p_{\B_j}
  -\smashoperator{\sum_{\substack{2\le i,j\le m \\ i+j=m+2}}} p_{\A_i\B_j}
  +\smashoperator{\sum_{\substack{1\le i,j\le m \\ i+j\le m+1}}} p_{\A_i\B_j}
  \le0,
  \label{eq:imm22}
\end{equation}
are valid for $\CorP(\K_{m,m})$ for all $m\ge1$, generalizing CHSH
inequality~\cite{ClaHorShiHol-PRL69} for $m=2$ which is a facet of
$\CorP(\K_{2,2})$.
They conjectured that for any $m\ge1$, the $\I_{mm22}$ inequality is a
facet of $\CorP(\K_{m,m})$, and showed that the conjecture is true for
$m\le7$.

\paragraph*{Triangular elimination}
Avis, Imai, Ito and Sasaki~\cite{AviImaItoSas:0404014} proposed
\emph{triangular elimination} operation to convert any facet inequality
of $\CutP_n$ other than the triangle inequality to a facet of
$\CutP(\K_{1,m,m})$ for appropriate $m$.
A basic step in this conversion is described in the following theorem.

\begin{thm}[\cite{AviImaItoSas:0404014}] \label{thm:trielim}
  Let $G=(V,E)$ be a graph and $uu'\in E$ an edge of $G$.
  Let $W\subseteq\N_G(u)\cap\N_G(u')$ be a set of nodes that are
  adjacent to both $u$ and $u'$.
  We define a graph $G^+=(V^+,E^+)$, the \emph{detour extension} of
  $G$, as follows.
  We add a new node $v$ to $G$ in the middle of the edge $uu'$,
  dividing $uu'$ into two edges $uv$ and $u'v$, and add new edges $vw$
  for each $w\in W$.

  Let $\vct{a}^\trans\vct{x}\le a_0$ be a facet inequality of
  $\CutP(G)$.
  Define $\vct{b}^\trans\vct{x}\le a_0$, the \emph{triangular
  elimination} of $\vct{a}^\trans\vct{x}\le a_0$, to be the inequality
  obtained by combining the triangle inequality
  $-a_{uu'}x_{uu'}+a_{uu'}x_{uv}-\lvert a_{uu'}\rvert x_{u'v}\le0$ with
  $\vct{a}^\trans\vct{p}\le a_0$.

  If there exists an edge $e\in E\setminus(\{uu'\}\cup\{uw,u'w\mid
  w\in W\})$ such that $a_e\ne0$, then the inequality
  $\vct{b}^\trans\vct{x}\le a_0$ is a facet of $\CutP(G^+)$.
\end{thm}

\section{Inequality $I(G,H)$: A generalization of $\Gr_7$}
  \label{sect:ineq-1}

In this section, we define the inequality $I(G,H)$ valid
for the cut polytope, and give a necessary and sufficient condition
for $I(G,H)$ to be a facet.

First we define the inequality.
Let $n\ge1$ be an integer, and $G=(V,E)$ and $H=(V,F)$ be two graphs
with $n$ nodes.
We require that the edges of $H$ are node-disjoint.
Let $t=\lvert F\rvert$ and $k=n-t$, and we denote the connected component
decomposition of $H$ by $V=V_1\cup\dots\cup V_k$.
Note that the size of any connected component $V_i$ is one or
two.
Finally we require that $E$ contains exactly $\binom{k}{2}$ edges:
for each $1\le i<j\le k$ there is an edge $e_{ij}$
connecting a
node in $V_i$ and a node in $V_j$.
We consider the following inequality which we denote as $I(G,H)$:
\begin{equation}
  \smashoperator[r]{\sum_{uv\in E}} T(u,v;n+1)
  -\smashoperator{\sum_{uv\in F}} T(u,v;n+1)
  +2\smashoperator{\sum_{V_i=\{u\}}} x_{u,n+1} \le2.
  \label{eq:ineq-1}
\end{equation}

For example, $I(\K_2,\overline{\K}_2)$ is identical to the triangle
inequality and $I(\K_4,\overline{\K}_4)$ to the pure pentagonal inequality,
where $K_n$ is the complete graph on $n$ nodes, and $\overline{\K}_n$
is its complement.

It is sometimes convenient to relabel the nodes in $V$ so that $H$ is
in a restricted form.
For $k\ge1$ and $0\le t\le k$, let $H_{k,t}=(V,E)$ be a graph with node set
$V=\{1,\dots,k+t\}$ and edge set $E=\{(i,k+i)\mid 1\le i\le t\}$.
Then any graph $H$ with $n=k+t$ nodes and $t$ node-disjoint edges can
be relabelled to $H_{k,t}$, and therefore we can restrict $I(G,H)$ to
$I(G,H_{k,t})$ without loss of generality.

\begin{figure}
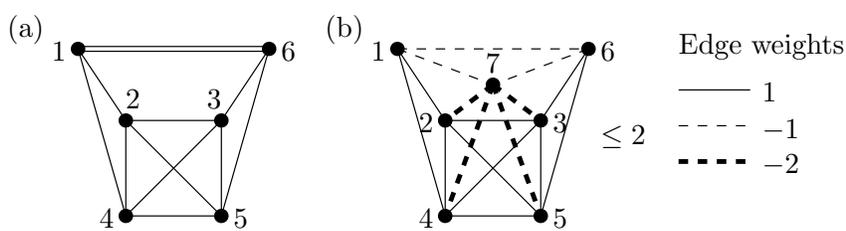

  \centering
  \begin{tabular}{ccc}
    (a) \raisebox{-\height}{\input{figures/g6.pstex_t}}
    &
    (b) \raisebox{-\height}{\input{figures/gr7-g-ineq.pstex_t}}%
        \hspace*{1em}
    &
    \raisebox{-\height}{\input{figures/edge-weights-1m1m2.pstex_t}}%
    \hspace*{1em}
  \end{tabular}
  \caption{(a)~A graph $G_6=(V,E)$ (edges drawn as single lines) and a
    graph $H_{5,1}=(V,F)$ (an edge drawn as a double line).
    (b)~The inequality $I(G_6,H_{5,1})$, which is a
    switching of the $\Gr_7$ inequality.}
  \label{fig:gr7-g}
\end{figure}

We check that the $\Gr_7$ inequality is a switching of an
inequality of this kind.
Let $G_6=(V,E)$ and $H_{5,1}=(V,F)$ be the graphs with six nodes shown in
Figure~\ref{fig:gr7-g}~(a).
Then the inequality $I(G_6,H_{5,1})$ is as shown in
Figure~\ref{fig:gr7-g}~(b).
We switch $I(G_6,H_{5,1})$ by the cut $\{1,6\}$ and change the
labels of nodes $1,2,3,4,5,6,7$ to $6,1,2,3,4,7,5$, respectively.
Then the resulting inequality is identical to $\Gr_7$.

Now we prove the validity of $I(G,H)$.

\begin{prop} \label{prop:valid-1}
  The inequality $I(G,H)$ is valid for $\CutP_{n+1}$.
  In addition, the cut vector $\vct{\delta}(S)$ with $S\subseteq V$ is a
  root of $I(G,H)$ if and only if one of the following
  conditions is satisfied.
  \begin{enumerate}[(i)]
    \item
      There exists a unique $i$ such that $V_i\subseteq S$, and no
      edge of $G$ is contained in $S$.
    \item
      There exist exactly two values of $i$ (let them be $i_1$ and
      $i_2$) such that $V_i\subseteq S$.
      In addition, $e_{i_1i_2}$ is the only edge of $G$ that is
      contained in $S$.
  \end{enumerate}
\end{prop}

\begin{proof}
  We show that the cut vector $\vct{\delta}(S)$ defined by any subset
  $S$ of $V$ satisfies the inequality $I(G,H)$.
  Note that with $\vct{x}=\vct{\delta}(S)$, each term evaluates to
  either to zero or two.

  Let $A=\{i\mid V_i\subseteq S\}$ and $B=\{ij\mid e_{ij}\subseteq
  S\}$.
  The left hand side of $I(G,H)$ evaluated with
  $\vct{x}=\vct{\delta}(S)$ is equal to $2\lvert A\rvert-2\lvert
  B\rvert$.
  Now $\lvert B\rvert\ge\binom{\lvert A\rvert}{2}$, since for each of
  the $\binom{\lvert A\rvert}{2}$ pairs $ij$ of elements of $A$, there
  is an edge $e_{ij}$ with both endpoints in $S$.
  Therefore we have $
        2\lvert A\rvert-2\lvert B\rvert
    \le 3\lvert A\rvert-\lvert A\rvert^2
    =2-(\lvert A\rvert-1)(\lvert A\rvert-2) \le 2
  $.
  So \eqref{eq:ineq-1} is valid.

  The condition for roots is obtained from the fact that this
  inequality is satisfied with equality if and only if $\lvert
  A\rvert$ is one or two and $\lvert B\rvert=\binom{\lvert
  A\rvert}{2}$.
\end{proof}

Now we consider when the inequality $I(G,H)$ becomes a facet of
$\CutP_{n+1}$.

\begin{thm} \label{thm:facet-1}
  Assume $k\ge3$.
  Then the inequality $I(G,H)$ is a facet of $\CutP_{n+1}$ if and only
  if all nodes in $G$ have degree at least two.
\end{thm}

\begin{proof}
  As mentioned above, we can assume $H=H_{k,t}$ without loss of
  generality.

  First we prove the ``only if'' part.
  Let $u$ be a node whose degree in $G$ is at most one.
  In this case $H_{k,t}$ has an edge incident to node $u$.
  Without loss of generality, we assume $u=k+t$.
  If the degree of node $k+t$ in $G$ is one, then let $v$ be
  the only node that is adjacent to node $k+t$ in $G$.
  Otherwise let $v=n+1$.
  In both cases, $I(G,H_{k,t})$ is the sum of a triangle inequality
  $T(k+t,v;t)\le0$ and the inequality $I(G/(t,k+t),H_{k,t-1})$, where
  $G/(t,k+t)$ is a graph obtained from $G$ by identifying two nodes
  $t$ and $k+t$ into a node $t$.
  Therefore, $I(G,H_{k,t})$ is not a facet of $\CutP_{n+1}$.

  Now we prove the ``if'' part.
  The proof is by induction on $t$.

  First we consider the case $t=0$.
  In this case, $H_{k,0}$ has no edges and $G$ is the complete graph
  $\K_n$.
  Switching the inequality $I(\K_n,\overline{\K}_n)$ by the cut $\{1\}$ gives
  a hypermetric inequality defined by an integer vector $\vct{b}$ with
  $b_{n+1}=-(k-3)$, $b_1=-1$ and $b_2=\dots=b_n=1$.
  This hypermetric inequality is a facet of $\CutP_{n+1}$ by
  Theorem~\ref{thm:hyp}.

  Now we consider the case $t\ge1$.
  Note that contracting the edge $(t,t+k)$ in $H_{k,t}$ gives
  $H_{k,t-1}$.
  Key facts are that the inequality $I(G,H_{k,t})$ is obtained by
  lifting $I(G/(t,t+k),H_{k,t-1})$, and that $I(G/(t,t+k),H_{k,t-1})$
  is a facet of $\CutP_n$ by the induction hypothesis.

  Let $V_i=V_{i+k}=\{i,i+k\}$
  for $1\le i\le t$ and $V_i=\{i\}$ for $t+1\le i\le k$.
  We define $n$ subsets of $V$ as follows.

  \begin{itemize}
    \setlength{\itemsep}{0pt}
    \item
      Let $p$ and $p'$ be two distinct nodes adjacent to node $t+k$ in
      $G$.
      Then define $T^{(1)}=\{t\}\cup V_p\cup V_{p'}$.
    \item
      Let $q$ and $q'$ be two distinct nodes adjacent to node $t$ in
      $G$.
      Then define $T^{(2)}=\{t+k\}\cup V_q\cup V_{q'}$.
    \item
      For each $1\le i\le k$ with $i\ne t$, we define $T^{(3)}_i$.
      If $e_{it}$ has an endpoint $t+k$, then $T^{(3)}_i=\{t\}\cup V_i$.
      Otherwise, $T^{(3)}_i=\{t+k\}\cup V_i$.
    \item
      For each $1\le i\le t-1$, we define a subset $T^{(4)}_i$.
      Let $u$ be either $i$ or $i+k$ that is an endpoint of the edge
      $e_{it}$, and $\bar{u}$ be either $i$ or $i+k$ that is different
      from $u$.
      Let $v$ be any node in $\N_G(u)\setminus V_t$ and choose $j$ so
      that $V_j\ni v$.
      Let $\bar{w}$ be
      either $t$ or $t+k$ that is not an endpoint of the edge
      $e_{jt}$.
      Then define $T^{(4)}_i=\{\bar{u},\bar{w}\}\cup V_j$.
  \end{itemize}

  It is easy to check that each of these subsets is a root of
  $I(G,H_{k,t})$ and contains exactly one of $t$ and $t+k$.
  Note that none of them contains node $n+1$.

  The following claim can be proved in a straightforward way.
  A proof is included in Appendix~\ref{apdx:incidence-independent-1}.

  \begin{cl} \label{cl:incidence-independent-1}
    The $n$ incident vectors of $T^{(1)}$, $T^{(2)}$,
    $T^{(3)}_i\;(i\ne t)$ and $T^{(4)}_i\;(1\le i\le t-1)$ are
    linearly independent.
  \end{cl}

  From now on, we refer to the $n$ sets $T^{(1)}$, $T^{(2)}$,
  $T^{(3)}_i\;(i\ne t)$ and $T^{(4)}_i\;(1\le i\le t-1)$ as
  $T_1,\dots,T_n$.

  Let $\vct{a}^\trans\vct{x}\le0$ be the switching of $I(G,H_{k,t})$ by
  its root $\{t,t+k\}$.

  The $(t,t+k)$-collapse $(\vct{a}^{t,t+k})^\trans\vct{x}\le0$ of
  $\vct{a}^\trans\vct{x}\le0$ is the switching by the cut $\{t\}$ of
  $I(G/(t,t+k),H_{k,t-1})$, which is a facet of $\CutP_n$ by induction
  hypothesis.

  For $1\le i\le n$, let $T'_i$ be $T_i\symdiff\{t,t+k\}$ if $t\in
  T_i$, and $(V\cup\{n+1\})\setminus(T_i\symdiff\{t,t+k\})$
  otherwise, where $\symdiff$ means the symmetric difference of two
  sets.
  Then $T'_i$ is a root of the inequality $\vct{a}^\trans\vct{x}\le0$
  and contains $t+k$ but does not contain $t$.
  In addition, the $n$ vectors $T'_1,\dots,T'_n$ are also
  linearly independent.
  From Lemma~\ref{lem:lifting}, the inequality
  $\vct{a}^\trans\vct{x}\le0$ is a facet of $\CutP_{n+1}$, which means
  $I(G,H_{k,t})$ is also a facet of $\CutP_{n+1}$.
\end{proof}

\begin{figure}
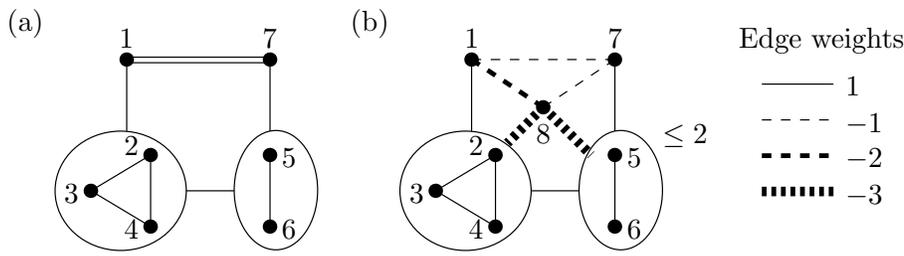

  \centering
  \begin{tabular}{ccc}
    (a)
    \raisebox{-\height}{\input{figures/facet-1-cut8-graph.pstex_t}}
    &
    (b)
    \raisebox{-\height}{\input{figures/facet-1-cut8-ineq.pstex_t}}%
    \hspace*{1.5em}
    &
    \hspace*{0.5em}%
    \raisebox{-\height}{\input{figures/edge-weights-1m1m2m3.pstex_t}}%
    \hspace*{0.5em}
  \end{tabular}
  \caption{(a)~A graph $G$ (edges drawn as single lines) and $H_{6,1}$ (an
    edge drawn as a double line).
    (b)~The inequality $I(G,H_{6,1})$, which is proved to be a
    facet of $\CutP_8$ by Theorem~\ref{thm:facet-1}.
    A line connected to a circle enclosing nodes $1$, $2$ and $3$
    represents $3$ edges with identical weights each connected to the
    nodes $1$, $2$ and $3$.
    Similar for lines connected to the other circles.}
  \label{fig:facet-1-cut8}
\end{figure}

For example, let us consider the graphs $G=(V,E)$ and $H_{6,1}=(V,F)$ shown
in Figure~\ref{fig:facet-1-cut8}~(a).
In this case the inequality $I(G,H_{6,1})$, illustrated in
Figure~\ref{fig:facet-1-cut8}~(b), is a facet of
$\CutP_8$ by Theorem~\ref{thm:facet-1}.

\section{Inequality $I'(G,H,C)$: A generalization of $\Gr_8$}
  \label{sect:ineq-2}

Let $G=(V,E)$, $H=(V,F)$, $n=\lvert V\rvert$, $t=\lvert F\rvert$,
$k=n-t$ and
$V=V_1\cup\dots\cup V_k$ be as defined in Section~\ref{sect:ineq-1}.
In this section we require an additional condition that $G$ has a
cycle $C$ of length four (this condition implies $k\ge4$).
Let $V_C$ be the set of the four nodes of $C$.
Then we consider an inequality for the cut polytope on $n+2$ nodes:
\begin{equation}
  \smashoperator[r]{\sum_{uv\in E}} T(u,v;n+1)
  -\smashoperator{\sum_{uv\in F}} T(u,v;n+1)
  +2\smashoperator{\sum_{V_i=\{u\}}} x_{u,n+1}
  +\smashoperator{\sum_{u\in V_C}} (x_{u,n+1}-x_{u,n+2}) \le2.
  \label{eq:ineq-2}
\end{equation}
We refer to inequality~\eqref{eq:ineq-2} by $I'(G,H,C)$.
Note that the $(n+1,n+2)$-collapsing of $I'(G,H,C)$ is identical to
$I(G,H)$.

\begin{figure}
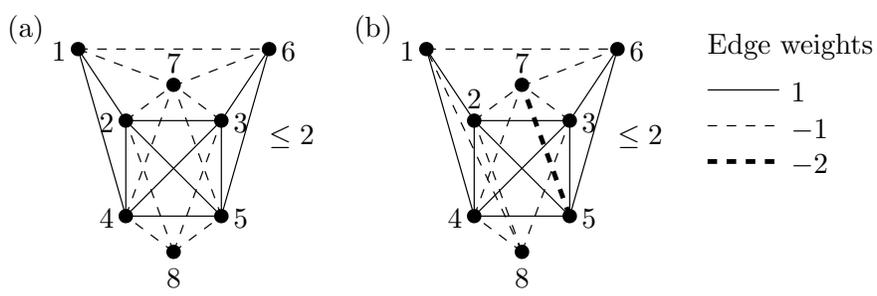

  \centering
  \begin{tabular}{ccc}
    (a) \raisebox{-\height}{\input{figures/gr8-g-ineq.pstex_t}}%
        \hspace*{1em}
    &
    (b) \raisebox{-\height}{\input{figures/facet-2-cut8-ineq.pstex_t}}%
        \hspace*{1em}
    &
    \raisebox{-\height}{\input{figures/edge-weights-1m1m2.pstex_t}}
  \end{tabular}
  \caption{Two inequalities $I'(G_6,H_{5,1},C)$ with different $C$.
    Both are proved to be facets of $\CutP_8$ by
    Theorem~\ref{thm:facet-2}.
    (a)~Case of $C=\{23,34,45,52\}$.
    The inequality is a switching of the $\Gr_8$ inequality.
    (b)~Case of $C=\{12,23,34,41\}$.}
  \label{fig:facet-2}
\end{figure}

As an example, we show that the $\Gr_8$ inequality is a switching of
an inequality of this kind.
Consider again the graphs $G_6=(V,E)$ and $H_{5,1}=(V,F)$ shown in
Figure~\ref{fig:gr7-g}~(a).
Note that $G_6$ contains a cycle $C=\{23,34,45,52\}$ of length four.
Then the inequality $I'(G,H,C)$ is as shown in
Figure~\ref{fig:facet-2}~(a), and switching it by the cut $\{1,6\}$ and
relabelling nodes appropriately gives the $\Gr_8$ inequality.

\begin{prop} \label{prop:valid-2}
  The inequality $I'(G,H,C)$ is valid for $\CutP_{n+2}$.
\end{prop}

\begin{proof}
  Let $M$ be a set of two node-disjoint edges in the cycle $C$.
  Note that there are two choices of $M$.
  No matter which set we choose as $M$, the inequality $I'(G,H,C)$ can be
  written as
  \begin{equation}
    \smashoperator[r]{\sum_{uv\in E\setminus M}}
      T(u,v;n+1)
    +\smashoperator{\sum_{uv\in M}}
      T(u,v;n+2)
    -\smashoperator{\sum_{uv\in F}}
      T(u,v;n+1)
    +2\smashoperator{\sum_{V_i=\{u\}}}
      x_{u,n+1}
    \le2.
    \label{eq:ineq-2-matching}
  \end{equation}

  We show that the cut vector $\vct{\delta}(S)$ defined by any subset
  $S$ of $V\cup\{n+2\}$ satisfies \eqref{eq:ineq-2-matching}.
  Let $A=\{i\mid V_i\subseteq S\}$ and
  $B=\{ij\mid e_{ij}\subseteq S,\; e_{ij}\in E\setminus M\}$.
  Now $\lvert B\rvert\ge\binom{\lvert A\rvert}{2}
  -\lfloor\lvert A\rvert/2\rfloor$, since for each
  $ij\in B$ there is an edge $e_{ij}$ with both endpoints in $A$,
  except for up to $\lfloor\lvert A\rvert/2\rfloor$
  edges that may be part of $M$.
  The left hand side of \eqref{eq:ineq-2} evaluated with
  $\vct{x}=\vct{\delta}(S)$ is at most $2\lvert A\rvert-2\lvert
  B\rvert$.
  Combining inequalities we have
  $
        2\lvert A\rvert-2\lvert B\rvert
    \le 3\lvert A\rvert+2\lfloor\lvert A\rvert/2\rfloor
        -\lvert A\rvert^2
    \le2
  $
  except when $\lvert A\rvert=2$.
  So \eqref{eq:ineq-2-matching} is valid for all these cases.
  Suppose $\lvert A\rvert=2$.

  \medskip
  \noindent
  Case 1: The two nodes in $A$ do not form an edge in $M$.

  \noindent
  In this case $\lvert B\rvert=1$, the LHS of
  \eqref{eq:ineq-2-matching} is at most $2\lvert A\rvert-2\lvert
  B\rvert=2$ and the inequality is valid.

  \medskip
  \noindent
  Case 2: The two nodes in $A$ form an edge in $M$.

  \noindent
  In this case we replace $M$ by $C\setminus M$.
  This does not change the LHS of \eqref{eq:ineq-2-matching}, and the
  inequality is valid by Case~1.
\end{proof}

Before we state a sufficient condition for $I'(G,H,C)$ to be a facet
of $\CutP_{n+2}$, we assume some conditions on $H$ and $C$ without
loss of generality.
We assume $H=H_{k,t}$, where $H_{k,t}$ is the same as that defined in
the previous section, and we also assume that indices of the four
nodes of $C$ are at most $k$.
We say that node $i$ in $C$ is \emph{free} if $1\le i\le
t$ and $i+k$ is incident to edge $e_{ij}$ where $j$ is the unique node
in $C$ that is not adjacent to $i$ in $C$.
The following theorem gives a sufficient condition for
$I'(G,H_{k,t},C)$ to be a facet.

\begin{thm} \label{thm:facet-2}
  The inequality $I'(G,H_{k,t},C)$ is a facet of $\CutP_{n+2}$ if all
  of the following conditions are satisfied:
  \begin{enumerate}[(i)]
    \item
      All nodes in $G$ have at least two neighbors.
    \item \label{enum:facet-2-condition-2}
      For each $t+1\le i\le k$ except for nodes in $C$, there exists a
      free node $j$ in $C$ such that $e_{ij}$ is incident to $j+k$.
    \item \label{enum:facet-2-condition-3}
      For each $1\le i\le t$ except for nodes in $C$, either:
      \begin{itemize}
        \item
          Nodes $i$ and $i+k$ are incident to exactly two out of four
          edges $e_{ij}$ with $j\in V_C$, or
        \item
          There exists a free node $j$ in $C$ such that $e_{ij}$ is
          incident to $j+k$.
      \end{itemize}
  \end{enumerate}
\end{thm}

Since $I'(G,H,C)$ is a lifting of $I(G,H)$, we may prove
Theorem~\ref{thm:facet-2} by combining the lifting lemma
(Lemma~\ref{lem:lifting}) with Theorem~\ref{thm:facet-1}.
The proof is given in Appendix~\ref{apdx:facet-2}.
As an example of the theorem, consider the graphs $G_6$ and $H_{5,1}$
shown in Figure~\ref{fig:gr7-g}~(a), but this time let
$C=\{12,23,34,41\}$.
In this case the inequality $I'(G,H,C)$, shown in
Figure~\ref{fig:facet-2}~(b), is a facet of $\CutP_8$ by
Theorem~\ref{thm:facet-2}.

Unlike $I(\K_k,\overline{\K}_k)$, which is always a facet of
$\CutP_{k+1}$, the face of $\CutP_{k+2}$ supported by the inequality
$I'(\K_k,\overline{\K}_k,C)$ with $k\ge5$ and $C=\{12,23,34,41\}$ is
contained in a triangle facet $x_{5,k+2}-x_{5,k+1}-x_{k+1,k+2}\le0$
and never supports a facet.

\section{Tightness of the $\I_{mm22}$ Bell inequalities}
  \label{sect:imm22}

In this section, we prove that for any $m$, the $\I_{mm22}$ inequality
is a facet of $\CorP(\K_{m,m})$, or in other words, a tight Bell
inequality.
Since the proof does not depend on the proof of validity given
in~\cite{ColGis-JPA04}, our proof also serves as another way to prove
the validity of the $\I_{mm22}$ inequality.

Let $\K_{1,m,m}=(V_{1,m,m},E_{1,m,m})$ be a complete tripartite graph
with node set $V_{1,m,m}=\{\Z,\A_1,\dots,\A_m,\B_1,\dots,\B_m\}$ and
edge set $E_{1,m,m}=\{\Z\A_i\mid 1\le i\le m\}\cup
\{\Z\B_j\mid 1\le j\le m\}\cup\{\A_i\B_j\mid 1\le i,j\le m\}$.
We rewrite the $\I_{mm22}$ inequality to an inequality for
$\CutP(\K_{1,m,m})$ by using the covariance mapping.
We switch this inequality by the cut $\{\A_1,\dots,\A_m\}$.
After that, we change the labels of the $m$ nodes
$\B_1,\B_2,\dots,\B_m$ to $\B_{m+1},\B_m,\dots,\B_2$, respectively,
both in the inequality and the graph $\K_{1,m,m}$.
Let us denote the resulting complete tripartite graph by
$\K'_{1,m,m}$.
Then the inequality becomes
\begin{equation}
  -(m-2) x_{\Z\A_1}
  -\smashoperator{\sum_{2\le i\le m}} (m-i) x_{\Z\A_i}
  -(m-2) x_{\Z\B_{m+1}}
  -\smashoperator{\sum_{2\le j\le m}} (j-2) x_{\Z\B_j}
  -\smashoperator{\sum_{2\le i\le m}} x_{\A_i\B_i}
  +\smashoperator{\sum_{1\le i<j\le m+1}} x_{\A_i\B_j}\le2.
  \label{eq:imm22-cut-2}
\end{equation}
It is easy to check that the inequality \eqref{eq:imm22-cut-2} is
identical to the inequality $I(G,H)$ with $G=(V,E)$ and $H=(V,F)$,
where $V=\{\A_1,\dots,\A_m,\B_2,\dots,\B_{m+1}\}$,
$E=\{\A_i\B_j\mid 1\le i<j\le m+1\}$, and
$F=\{\A_i\B_i\mid 2\le i\le m\}$.
Therefore the $\I_{mm22}$ inequality is a tight Bell inequality if
and only if the inequality $I(G,H)$ is a facet of
$\CutP(\K'_{1,m,m})$.

Note that we cannot use Theorem~\ref{thm:facet-1} directly to prove
that $I(G,H)$ is a facet, since the graph $G$ does not
satisfy the condition of Theorem~\ref{thm:facet-1}.
However, if we assume $m\ge3$, the inequality $I(G,H)$ is the
triangular elimination of another inequality $I(G',H')$, where $G'$
(resp.\ $H'$) is the graph obtained from $G$ (resp.\ $H$) by
identifying node $\B_2$ to $\A_2$ and $\A_m$ to $\B_m$.
The inequality $I(G',H')$ is proved to be a facet of $\CutP_{2m-1}$ by
Theorem~\ref{thm:facet-1}.
Now, as was pointed out in \cite{AviImaItoSas:0404014}
and \cite{ItoSasImaAvi-EQIS04},
we can apply triangular elimination twice to the facet inequality
$I(G',H')$ to obtain $I(G,H)$.
The first application is done with $uu'=\A_1\A_2$, $v=\B_2$ and
$W=\{\Z,\A_3,\dots,\A_{m-1},\B_m,\B_{m+1}\}$.
The second application is done with $uu'=\B_m\B_{m+1}$, $v=\A_m$
and $W=\{\Z,\B_2,\dots,\B_{m-1}\}$.
Therefore, from Theorem~\ref{thm:trielim}, $I(G,H)$ is a facet of
$\CutP(\K'_{1,m,m})$.

Since it is easy to check the cases $m=1$ and $2$, we obtain the
following theorem.

\begin{thm} \label{thm:imm22}
  For any $m\ge1$, the $\I_{mm22}$ inequality~\eqref{eq:imm22} is a
  tight Bell inequality.
\end{thm}

\section*{Acknowledgements}

We thank Hiroshi Imai and Yuuya Sasaki for useful discussions.

\bibliography{bell}

\providecommand{\SU}{{\mathrm{SU}}}
\begin{thebibliography}{10}

\bibitem{AviDez-Net91}
\textsc{D.~Avis and M.~Deza},
\newblock The cut cone, {$\mathrm{L}^1$} embeddability, complexity and
  multicommodity flows,
\newblock {\em Networks} (1991) \textbf{21}:595--617.

\bibitem{AviImaItoSas:0404014}
\textsc{D.~Avis, H.~Imai, T.~Ito, and Y.~Sasaki},
\newblock Deriving tight {B}ell inequalities for 2 parties with many 2-valued
  observables from facets of cut polytopes,
\newblock arXiv:quant-ph/0404014,  (2004).

\bibitem{ClaHorShiHol-PRL69}
\textsc{J.~F. Clauser, M.~A. Horne, A.~Shimony, and R.~A. Holt},
\newblock Proposed experiment to test local hidden-variable theories,
\newblock {\em Physical Review Letters} (1969) \textbf{23}(15):880--884.

\bibitem{ColGis-JPA04}
\textsc{D.~Collins and N.~Gisin},
\newblock A relevant two qubit {B}ell inequality inequivalent to the {CHSH}
  inequality,
\newblock {\em Journal of Physics A: Mathematical and General} (2004)
  \textbf{37}(5):1775--1787,
\newblock arXiv:quant-ph/0306129.

\bibitem{DesDezLau-DM94}
\textsc{C.~{\uppercase{d}e}~Simone, M.~Deza, and M.~Laurent},
\newblock Collapsing and lifting for the cut cone,
\newblock {\em Discrete Mathematics} (1994) \textbf{127}(1--3):105--130.

\bibitem{DezLau-JCAM94:applications1}
\textsc{M.~Deza and M.~Laurent},
\newblock Applications of cut polyhedra {I},
\newblock {\em Journal of Computational and Applied Mathematics} (1994)
  \textbf{55}(2):191--216.

\bibitem{DezLau-JCAM94:applications2}
\textsc{M.~Deza and M.~Laurent},
\newblock Applications of cut polyhedra {II},
\newblock {\em Journal of Computational and Applied Mathematics} (1994)
  \textbf{55}(2):217--247.

\bibitem{DezLau:cut97}
\textsc{M.~M. Deza and M.~Laurent},
\newblock {\em Geometry of Cuts and Metrics}, volume~15 of {\em Algorithms and
  Combinatorics},
\newblock Springer,  (1997).

\bibitem{Gri-EJC90}
\textsc{V.~P. Grishukhin},
\newblock All facets of the cut cone {$C_n$} for {$n=7$} are known,
\newblock {\em European Journal of Combinatorics} (1990) \textbf{11}:115--117.

\bibitem{ItoSasImaAvi-EQIS04}
\textsc{T.~Ito, Y.~Sasaki, H.~Imai, and D.~Avis},
\newblock Families of tight {B}ell inequalities derived from classes of facets
  of cut polytopes,
\newblock In {\em Proceedings of ERATO conference on Quantum Information
  Science (EQIS'04)}, pages 78--79,  (2004).

\bibitem{WerWol-QIC01}
\textsc{R.~F. Werner and M.~M. Wolf},
\newblock Bell inequalities and entanglement,
\newblock {\em Quantum Information \& Computation} (2001) \textbf{1}(3):1--25,
\newblock arXiv:quant-ph/0107093.

\end{thebibliography}

\appendix

\section{Proof of Claim~\ref{cl:incidence-independent-1}}
  \label{apdx:incidence-independent-1}

\begin{proof}[Proof of Claim~\ref{cl:incidence-independent-1}]
  First let $1\le i\le t-1$.
  In these $n$ sets, $T^{(4)}_i$ is the only one that contains
  exactly one of $i$ and $i+k$.
  This means that the linear independence of the incident vectors of
  $k+1$ sets $T^{(1)}$, $T^{(2)}$ and $T^{(3)}_i\;(i\ne t)$ implies
  the linear independence of all the $n$ incident vectors.

  Next let $1\le i\le k$, $i\ne t,p,p',q,q'$.
  In these $k+1$ sets, $T^{(3)}_i$ is the only one that contains
  $i$.
  This means that the linear independence of the incident vectors of
  6 sets $T^{(1)}$, $T^{(2)}$, $T^{(3)}_p$, $T^{(3)}_{p'}$,
  $T^{(3)}_q$ and $T^{(3)}_{q'}$ implies the linear independence of
  all the $n$ incident vectors.

  Finally, these 6 incident vectors are linearly independent since
  they form an $(n+1)\times6$ matrix containing 6 rows which form a
  nonsingular matrix:
  \[
    \left(\begin{array}{ccccccc}
      (t)   & 1 & 1 & 1 & 0 & 0 & 0 \\
      (t+k) & 0 & 0 & 0 & 1 & 1 & 1 \\
      (p)   & 1 & 1 & 0 & 0 & 0 & 0 \\
      (p')  & 1 & 0 & 1 & 0 & 0 & 0 \\
      (q)   & 0 & 0 & 0 & 1 & 1 & 0 \\
      (q')  & 0 & 0 & 0 & 1 & 0 & 1
    \end{array}\right).  \qedhere
  \]
\end{proof}

\section{Proof of Theorem~\ref{thm:facet-2}}
  \label{apdx:facet-2}

First we note that from the proof of Proposition~\ref{prop:valid-2},
some of the roots of $I'(G,H,C)$ are characterized as follows.

\begin{prop} \label{prop:roots-2}
  A cut vector $\vct{\delta}(S)$ with $S\subseteq V\cup\{n+2\}$ is a
  root of $I'(G,H,C)$ if one of the following
  conditions is satisfied.
  \begin{enumerate}[(i)]
    \item
      $S$ does not contain node $n+2$, and $\vct{\delta}(S)$ is a
      root of $I(G,H)$.
    \item
      $S$ contains node $n+2$ and exactly two out of four nodes of $C$
      (possibly along with other nodes), and
      $\vct{\delta}(S\setminus\{n+2\})$ is a root of $I(G,H)$.
    \item
      $S=V_{c_1}\cup V_{c_2}\cup V_{c_3}\cup\{n+2\}$, where each
      $V_{c_i}$ $(i=1,2,3)$ contains a node of $C$.
  \end{enumerate}
\end{prop}

By using this characterization, we prove Theorem~\ref{thm:facet-2} as
follows.

\begin{proof}[Proof of Theorem~\ref{thm:facet-2}]
  The $(n+1,n+2)$-collapse of $I'(G,H_{k,t},C)$ is the inequality
  $I(G,H_{k,t})$, and it is a facet of
  $\CutP_{n+1}$ by Theorem~\ref{thm:facet-1}.

  Let $C=\{c_1c_2,c_2c_3,c_3c_4,c_4c_1\}$.
  Note that $1\le c_1,c_2,c_3,c_4\le k$.
  For $1\le\lambda,\mu\le4$, we denote by $\lambda\oplus\mu$ the
  unique integer $\nu$ such that $1\le\nu\le4$ and
  $\nu\equiv\lambda+\mu\pmod{4}$.
  Let $V_i=\{i,\,i+k\}$ for $1\le i\le t$ and $V_i=\{i\}$ for $t+1\le
  i\le k$.

  To use Lemma~\ref{lem:lifting}, we define $n+1$ subsets of
  $V\cup\{n+2\}$ as follows.
  \begin{itemize}
    \item
      We define $T^{(1)}_1=V_{c_1}\cup V_{c_2}\cup\{n+2\}$,
      $T^{(1)}_2=V_{c_1}\cup V_{c_3}\cup\{n+2\}$,
      $T^{(1)}_3=V_{c_1}\cup V_{c_4}\cup\{n+2\}$,
      $T^{(1)}_4=V_{c_2}\cup V_{c_3}\cup\{n+2\}$ and
      $T^{(1)}_5=V_{c_1}\cup V_{c_2}\cup V_{c_3}\cup\{n+2\}$.
    \item
      For each $1\le\lambda\le4$ such that $1\le c_\lambda\le t$, we
      define $T^{(2)}_\lambda=V_{c_{\lambda\oplus1}}\cup
      V_{c_{\lambda\oplus3}}\cup\{c_\lambda+k,\,n+2\}$.
    \item
      For each $t+1\le i\le k$ that is not incident to $C$, by
      condition~(\ref{enum:facet-2-condition-2}), there exists
      a free node $c_\lambda$ of $C$ such that $e_{ic_\lambda}$ is
      incident to $c_\lambda+k$.
      Then we define $T^{(3)}_i=V_{c_{\lambda\oplus2}}\cup
      \{c_\lambda,\,i,\,n+2\}$.
    \item
      For each $1\le i\le t$ that is not incident to $C$, we define
      $T^{(4)}_i$ and $T^{(5)}_i$
      as follows.
      \begin{itemize}
        \item
          If nodes $i$ and $i+k$ are incident to exactly two out of
          four edges $e_{ic}$ with $c\in V_C$, then define
          $\{\lambda,\lambda',\mu,\mu'\}=\{1,2,3,4\}$ such that
          $e_{ic_\lambda}$ and $e_{ic_{\lambda'}}$ are incident to $i$
          and $e_{ic_\mu}$ and $e_{ic_{\mu'}}$ are incident to $i+k$.
          In this case, we define
          $T^{(4)}_i=\{i,\,n+2\}\cup V_{c_\mu}\cup V_{c_{\mu'}}$ and
          $T^{(5)}_i=\{i+k,\,n+2\}\cup V_{c_\lambda}\cup
          V_{c_{\lambda'}}$.
        \item
          If not, let $u$ be either $i$ or $i+k$ that is incident to
          at most one of $e_{ic}$ with $c\in V_C$.
          By condition~(\ref{enum:facet-2-condition-3}), there exists
          a free node $c_\lambda$ such that $e_{ic_\lambda}$ is
          incident to $c_\lambda+k$.
          Then we define
          $T^{(4)}_i=V_{c_{\lambda\oplus2}}\cup
          \{c_\lambda,\,i,\,i+k,\,n+2\}$ and
          $T^{(5)}_i=V_\mu\cup V_{\mu+2}\cup\{u,\,n+2\}$, where $\mu$ is
          either $1$ or $2$ such that neither $e_{ic_\mu}$ nor
          $e_{ic_{\mu+2}}$ is incident to $u$.
      \end{itemize}
  \end{itemize}

  Each of these subsets contains $n+2$ but not $n+1$.
  By using Propositions~\ref{prop:valid-1} and \ref{prop:roots-2}, it
  is easy to check they are roots of $I'(G,H_{k,t},C)$.

  Now we prove the following claim.

  \begin{cl} \label{cl:incidence-independent-2}
    The $n+1$ incident vectors of $T^{(1)}_i$,
    $T^{(2)}_\lambda$, $T^{(3)}_i$,
    $T^{(4)}_i$ and $T^{(5)}_i$ are
    linearly independent.
  \end{cl}

  \begin{quote}
    \setlength{\parskip}{0pt}
    \setlength{\parindent}{\parindentorig}
    \begin{proof}
      The proof goes similarly to that of
      Claim~\ref{cl:incidence-independent-1}.

      For $t+1\le i\le k$ such that node $i$ is not incident to $C$,
      $T^{(3)}_i$ is the only set that includes $i$.
      Therefore all we have to prove is linear independence
      of the incidence vectors of sets $T^{(1)}_i$,
      $T^{(2)}_\lambda$, $T^{(4)}_i$ and $T^{(5)}_i$.

      For $1\le i\le t$ such that node $i$ is not incident to $C$,
      there are two possibilities.
      If nodes $i$ and $i+k$ are incident to exactly two out of four
      edges $e_{ci}$ with $c\in V_C$,
      the sets $T^{(4)}_i$ and $T^{(5)}_i$ are the only
      ones that include $i$ and $i+k$, respectively.
      Otherwise, the set $T^{(5)}$ is
      the only one that includes exactly one of two nodes $i$ and
      $i+k$, and $T^{(4)}_i$ is the only one that includes the other
      node in $i$ and $i+k$.
      Therefore, linear independence of the incidence vectors of
      sets $T^{(1)}_i$ and $T^{(2)}_\lambda$ will imply linear independence
      of all the $n+1$ incidence vectors.

      Let $1\le\lambda\le4$ such that $1\le c_\lambda\le t$.
      Among the remaining sets $T^{(1)}_i$ and $T^{(2)}_\mu$, the set
      $T^{(2)}_\lambda$ is the only one that includes exactly one of
      $c_\lambda$ and $c_\lambda+k$.

      Finally, the five incidence vectors of
      sets $T^{(1)}_1$, $T^{(1)}_2$, $T^{(1)}_3$, $T^{(1)}_4$ and
      $T^{(1)}_5$ are linearly independent since $(n+2)\times5$ matrix
      formed by these five vectors contains five rows with the pattern:
      \[
        \left(\begin{array}{cccccc}
          (c_1) & 1 & 1 & 1 & 0 & 1 \\
          (c_2) & 1 & 0 & 0 & 1 & 1 \\
          (c_3) & 0 & 1 & 0 & 1 & 1 \\
          (c_4) & 0 & 0 & 1 & 0 & 0 \\
          (n+2) & 1 & 1 & 1 & 1 & 1
        \end{array}\right),
      \]
      which is nonsingular.
    \end{proof}
  \end{quote}

  By Claim~\ref{cl:incidence-independent-2} and
  Lemma~\ref{lem:lifting}, $I'(G,H_{k,t},C)$ is a facet of $\CutP_{n+2}$.
\end{proof}

\end{document}